\documentclass[12pt]{article}
\usepackage{amsmath,amsfonts,amssymb,rotating,colordvi}
\usepackage{colordvi}
\def\qed{\nopagebreak\hfill{\rule{4pt}{7pt}}}
\def\proof{\noindent {\it{Proof.} \hskip 2pt}}

\newtheorem{theo}{Theorem}[section]

\newtheorem{lemm}[theo]{Lemma}

\newdimen\Squaresize \Squaresize=11pt
\newdimen\Thickness \Thickness=0.7pt
\def\Square#1{\hbox{\vrule width \Thickness
   \vbox to \Squaresize{\hrule height \Thickness\vss
    \hbox to \Squaresize{\hss#1\hss}
   \vss\hrule height\Thickness}
\unskip\vrule width \Thickness} \kern-\Thickness}

\def\Vsquare#1{\vbox{\Square{$#1$}}\kern-\Thickness}

\def\moins{\raise 1pt\hbox{{$\scriptstyle -$}}}

\begin{document}

\parskip 8pt

\begin{center}
{\large \bf On the Modes of Polynomials

             Derived from Nondecreasing Sequences}
\end{center}

\begin{center}
Donna Q. J. Dou$^{1}$, Arthur
L. B. Yang$^{2}$\\[6pt]

$^{1}$School of Mathematics\\
Jilin University, Changchun 130012, P. R. China

$^{2}$Center for Combinatorics, LPMC-TJKLC\\
Nankai University, Tianjin 300071, P. R. China

 Email: $^{1}${\tt qjdou@jlu.edu.cn}, $^{2}${\tt yang@nankai.edu.cn}
\end{center}

\noindent\textbf{Abstract.} Wang and Yeh proved that if $P(x)$ is a polynomial with nonnegative and nondecreasing coefficients, then $P(x+d)$ is unimodal for any $d>0$. A mode of a unimodal polynomial
$f(x)=a_0+a_1x+\cdots + a_mx^m$ is an index $k$
such that $a_k$ is the maximum coefficient.
Suppose that $M_*(P,d)$ is the smallest mode of $P(x+d)$, and $M^*(P,d)$ the greatest mode. Wang and Yeh conjectured that if $d_2>d_1>0$, then $M_*(P,d_1)\geq M_*(P,d_2)$ and $M^*(P,d_1)\geq M^*(P,d_2)$. We give a proof of this conjecture.



\noindent {\bf Suggested Running Title:} The mode conjecture

\section{Introduction}

This paper is concerned with the modes of unimodal polynomials constructed from nonnegative and nondecreasing sequences. Recall that a sequence $\{a_i\}_{0\leq i\leq m}$ is unimodal if there exists an index $0\leq k\leq m$ such that
$$a_0\leq \cdots\leq a_{k-1}\leq a_k\geq a_{k+1}\geq \cdots \geq a_m.$$
Such an index $k$ is called a mode of the sequence. Note that a mode of a sequence may not be unique. It is said to be
 {spiral} if
\begin{equation}
 a_m\leq a_0\leq a_{m-1}\leq a_1 \leq \cdots \leq a_{[\frac{m}{2}]},
\end{equation}
where $[\frac{m}{2}]$ stands for the greatest integer less than
$\frac{m}{2}$. Clearly, the spiral property implies unimodality. We say that a sequence  $\{a_i\}_{0\leq i\leq m}$ is {log-concave} if
for $1\leq k\leq m-1$,
$$a_{k}^2\geq a_{k+1}a_{k-1},$$
and it is {ratio monotone} if
\begin{equation}
\frac{a_m}{a_0} \leq \frac{a_{m-1}}{a_1} \leq \cdots\leq
 \frac{a_{m-i}}{a_i} \leq \cdots \leq \frac{a_{m-[\frac{m-1}{2}]}}{a_{[\frac{m-1}{2}]}} \leq 1
\end{equation}
and
\begin{equation}
\frac{a_0}{a_{m-1}} \leq \frac{a_{1}}{a_{m-2}} \leq \cdots\leq
 \frac{a_{i-1}}{a_{m-i}} \leq \cdots \leq \frac{a_{[\frac{m}{2}]-1}}{a_{m-[\frac{m}{2}]}} \leq 1.
\end{equation}
It is easily checked that the {ratio monotonicity} implies both log-concavity and the spiral property.

Let
$P(x)=a_0+a_1x+\cdots +a_mx^m$
be a polynomial with nonnegative coefficients. We say that $P(x)$ is unimodal if the sequence $\{a_i\}_{0\leq i\leq m}$ is unimodal. A mode of $\{a_i\}_{0\leq i\leq m}$ is also called a mode of $P(x)$. Similarly, we say that $P(x)$ is log-concave or ratio monotone if the sequence $\{a_i\}_{0\leq i\leq m}$ is log-concave or ratio monotone.

Throughout this paper $P(x)$ is assumed to be a polynomial with nonnegative and nondecreasing coefficients. Boros and Moll \cite{bormol1999} proved that $P(x + 1)$, as a polynomial of $x$, is unimodal. Alvarez et al. \cite{aabkmr2001} showed that
$P(x + n)$ is also unimodal for any positive integer $n$, and conjectured that $P(x+d)$ is unimodal for any $d>0$.
Wang and Yeh \cite{wangye2005} confirmed this conjecture and studied the modes of $P(x+d)$.
 Llamas and Mart\'{\i}nez-Bernal \cite{llamab2010} obtained the log-concavity of  $P(x + c)$ for
 $c\geq 1$. Chen, Yang and Zhou \cite{chenyangzhou} showed that $P(x+1)$ is ratio monotone, which leads to an alternative proof of the ratio monotonicity of the Boros-Moll polynomials \cite{chenxia08}.

Let $M_*(P,d)$ and $M^*(P,d)$ denote the smallest and the greatest mode of $P(x+d)$ respectively.
Our main result is the following theorem, which was conjectured by Wang and Yeh \cite{wangye2005}.

\begin{theo}\label{mainconj}
Suppose that $P(x)$
 is a monic polynomial of degree $m\geq 1$ with nonnegative and nondecreasing coefficients. Then for $0<d_1<d_2$, we have $M_*(P,d_1)\geq M_*(P,d_2)$ and $M^*(P,d_1)\geq M^*(P,d_2)$.
\end{theo}

From now on, we further assume that $P(x)$ is
monic, that is $a_m=1$.
For $0\leq k\leq m$, let
\begin{equation}\label{eqbk}
b_k(x)=\sum_{j=k}^m\binom{j}{k}a_jx^{j-k}.
\end{equation}
Therefore, $b_k(x)$ is of degree $m-k$ and $b_k(0)=a_k$.
For $1\leq k\leq m$, let
\begin{equation}\label{eqfk}
f_k(x)=b_{k-1}(x)-b_k(x),
\end{equation}
which is of degree $m-k+1$. Let $f_k^{(n)}(x)$ denote the $n$-th derivative of $f_k(x)$.

Our proof of Theorem \ref{mainconj} relies on
the fact that  $f_k(x)$ has only one real zero on $(0,+\infty)$. In fact, the derivative $f_k^{(n)}(x)$ of order $n\leq m-k$ has the same property. We establish this property by induction on $n$.

\section{Proof of Theorem \ref{mainconj}}


To prove Theorem \ref{mainconj}, we need the following three lemmas.

\begin{lemm} For any $0\leq k\leq m$, we have
$b_k'(x)=(k+1)b_{k+1}(x)$.
\end{lemm}

\proof It can be checked that
$$
\begin{array}{rcl}
b_k'(x) & = & \sum\limits_{j=k}^m\binom{j}{k}a_j(x^{j-k})'\\[8pt]
        & = & \sum\limits_{j=k+1}^m(j-k)\binom{j}{k}a_jx^{j-k-1}\\[8pt]
        & = & \sum\limits_{j=k+1}^m(j-k)\frac{j!}{k!(j-k)!}a_jx^{j-(k+1)}\\[8pt]
        & = & \sum\limits_{j=k+1}^m\frac{j!}{k!(j-k-1)!}a_jx^{j-(k+1)}\\[8pt]
        & = & \sum\limits_{j=k+1}^m(k+1)\frac{j!}{(k+1)!(j-(k+1))!}a_jx^{j-(k+1)}\\[8pt]
        & = & (k+1)b_{k+1}(x),
\end{array}
$$
as required. \qed

\begin{lemm}\label{mainlem} For $n\geq 1$ and $1\leq k\leq m$, we have
\begin{equation}\label{eq-nd}
f_k^{(n)}(x)=(k+n-1)_{n}b_{k+n-1}(x)-(k+n)_nb_{k+n}(x),
\end{equation}
where $(m)_j=m(m-1)\cdots (m-j+1)$.
\end{lemm}

\proof  Use induction on $n$. For $n=1$, we have
$$f_k^{(n)}(x)=f'(x)=kb_k-(k+1)b_{k+1}.$$
Assume that the lemma holds for $n=j$, namely,
$$f_k^{(j)}(x)=(k+j-1)_{j}b_{k+j-1}(x)-(k+j)_jb_{k+j}(x).$$
Therefore,
$$
\begin{array}{rcl}
f_k^{(j+1)}(x) & = & (k+j-1)_{j}b_{k+j-1}'(x)-(k+j)_jb_{k+j}'(x)\\[8pt]
             & = & (k+j)(k+j-1)_{j}b_{k+j}(x)-(k+j+1)(k+j)_jb_{k+j+1}(x)\\[8pt]
             & = & (k+j)_{j+1}b_{k+j}(x)-(k+j+1)_{j+1}b_{k+j+1}(x).
\end{array}
$$
This completes the proof.
\qed

\begin{lemm} For $1\leq k\leq m$ and $0\leq n\leq m-k$, the polynomial $f_k^{(n)}(x)$ has only one real zero on
the interval $(0,+\infty)$. In particular, $f_k(x)$ has only one real zero on
the interval $(0,+\infty)$.
\end{lemm}

\proof Use induction on $n$ from $m-k$ to $0$.
First, we consider the case $n=m-k$.
Recall that
$$f_k(x)=\sum_{j=k-1}^m \binom{j}{k-1}a_jx^{j-k+1}
-\sum_{j=k}^m\binom{j}{k}a_jx^{j-k}.$$
Thus  $f_k(x)$ is a polynomial of degree $m-k+1$.
Note that
$$f_k^{(m-k)}(x)=(m-k+1)!\binom{m}{k-1}a_m x+
\left[\binom{m-1}{k-1}a_{m-1}-\binom{m}{k}a_{m}\right](m-k)!.$$
Clearly, $f_k^{(m-k)}(x)$ has only one real zero $x_0$ on $(0,+\infty)$. So the lemma is true for
$n=m-k$.

Suppose that the lemma holds for $n=j$, where $m-k\geq j\geq 1$. We proceed to show that $f_k^{(j-1)}(x)$ has only one real zero on  $(0,+\infty)$.
From the inductive hypothesis it follows that
$f_k^{(j)}(x)$ has only one real zero on $(0,+\infty)$.
In light of \eqref{eq-nd}, it is easy to verify that $f_k^{(j)}(+\infty)>0$ and
$$f_k^{(j)}(0)= (k+j-1)_{j}a_{k+j-1}-(k+j)_ja_{k+j}\leq 0.$$
It follows that the polynomial $f_k^{(j-1)}(x)$
is decreasing up to certain point and becomes
increasing on the interval $(0,+\infty)$.
Again by \eqref{eq-nd} we find $f_k^{(j-1)}(+\infty)>0$ and
$$f_k^{(j-1)}(0)= (k+j-2)_{j-1}a_{k+j-2}-(k+j-1)_{j-1}a_{k+j-1}\leq 0.$$
So we conclude that $f_k^{(j-1)}(x)$ has only one real zero on  $(0,+\infty)$. This completes the proof.
\qed

\noindent \textit{Proof of Theorem \ref{mainconj}.}
In view of \eqref{eqbk}, we have
$$P(x+d)=\sum_{k=0}^ma_k(x+d)^k=\sum_{k=0}^m b_k(d)x^k.$$

Let us first prove that $M^*(P,d_1)\geq M^*(P,d_2)$. Suppose that $M^*(P,d_1)=k$.
If $k=m$, then the inequality $M^*(P,d_1)\geq M^*(P,d_2)$ holds.
For the case $0\leq k<m$, it suffices to verify that $b_k(d_2)> b_{k+1}(d_2)$.
By Lemma \ref{mainlem}, $f_{k+1}(x)$ has only one real zero on $(0,+\infty)$. Note that
$$f_{k+1}(0)\leq 0 \quad \mbox{ and } \quad f_{k+1}(+\infty)>0.$$
From $M^*(P,d_1)=k$ it follows that $b_k(d_1)> b_{k+1}(d_1)$, that is $f_{k+1}(d_1)>0$.
Therefore, $f_{k+1}(d_2)>0$, that is, $b_k(d_2)> b_{k+1}(d_2)$.

Similarly, it can be seen  that $M_*(P,d_1)\geq M_*(P,d_2)$.
Suppose that $M_*(P,d_2)=k$. If $k=0$, then we have $M_*(P,d_1)\geq M_*(P,d_2)$.
If $0<k\leq m$, it is necessary to show that $b_{k-1}(d_1)<b_k(d_1)$.
Again, by Lemma \ref{mainlem}, we know that $f_{k}(x)$ has only one real zero on $(0,+\infty)$.
From $M_*(P,d_2)=k$, it follows that $b_{k-1}(d_2)<b_k(d_2)$, that is $f_{k}(d_2)<0$.
By the boundary conditions
$$f_{k}(0)\leq 0  \quad \mbox{ and } \quad f_{k}(+\infty)>0,$$
we obtain $f_{k}(d_1)<0$, that is $b_{k-1}(d_1)<b_k(d_1)$.
This completes the proof.\qed

\noindent {\bf Acknowledgments.} This work was supported by the 973
Project, the PCSIRT Project of the Ministry of Education, and the
National Science Foundation of China.

\end{document}